%% file: hhv.tex
\documentclass[10pt]{amsart}
\usepackage{epic,eepic,latexsym}

 \newlength{\baseunit}               
 \newcount{\numlines}                
\newcommand{\point}{\vspace{3mm}\par\refstepcounter{subsection}{\bf \thesubsection.} }
\newcommand{\tpoint}[1]{\vspace{3mm}\par\refstepcounter{subsection}{\bf \thesubsection.} 
  {\em #1. ---} }
\newcommand{\epoint}[1]{\vspace{3mm}\par\refstepcounter{subsection}{\bf \thesubsection.} 
  {\em #1.} }
\newcommand{\bpoint}[1]{\vspace{3mm}\par\refstepcounter{subsection}{\bf \thesubsection.} 
  {\bf #1.} }

\newcommand{\epf}{\qed \vspace{+10pt}}

\newcommand{\E}{\mathbb{E}}

\newcommand{\com}{\mathbb{C}}

\newcommand{\proj}{\mathbb P}
\newcommand{\oh}{{\mathcal{O}}}

\newcommand{\cm}{{\mathcal{M}}}
\newcommand{\cmbar}{{\overline{\cm}}}

\newcommand{\al}{\alpha}

\newcommand{\Ga}{\Gamma}

\newcommand{\si}{\sigma}
\newcommand{\la}{\lambda}
\newcommand{\vir}{\operatorname{vir}}
\newcommand{\rbr}{\operatorname{br}}
\newcommand{\bbr}{\operatorname{Br}}

\newcommand{\Sym}{\operatorname{Sym}}
\newcommand{\Aut}{\operatorname{Aut}}

\newcommand{\etale}{\'{e}tale }
 
\newcommand{\cited}{}

\newcommand{\remind}[1]{{}}

\begin{document}
\pagestyle{plain}
\title{Hodge integrals and hurwitz numbers via virtual localization
\footnote{1991 Mathematics Subject Classification:  Primary 14H10, Secondary 14H30, 58D29}}
\author{Tom Graber}
\address{Dept. of Mathematics, Harvard University, Cambridge MA~02138} 
\email{graber@math.harvard.edu}
\author{Ravi Vakil}
\address{Dept. of Mathematics, MIT, Cambridge MA~02139}
\email{vakil@math.mit.edu}
\thanks{The second author is partially supported by NSF Grant DMS--9970101}
\date{February 29, 2000.}
\begin{abstract}
Ekedahl, Lando, Shapiro, and Vainshtein announced a
remarkable formula (\cite{elsv}) expressing Hurwitz numbers (counting covers of
$\proj^1$ with specified simple branch points, and specified branching
over one other point) in terms of Hodge integrals.  We give a proof of
this formula using virtual localization on the moduli space of stable
maps, and describe how the proof could be simplified by the proper
algebro-geometric definition of a ``relative space''.
\end{abstract}
\maketitle
\tableofcontents

{\parskip=12pt 

\section{Introduction}

Hurwitz numbers count certain covers of the projective line (or,
equivalently, factorizations of permutations into transposition).  They
have been studied extensively since the time of Hurwitz, and have
recently been the subject of renewed interest in physics (\cite{ct}),
combinatorics (\cite{d}, \cite{a}, and the series starting with
\cite{gj0}), algebraic geometry (recursions from Gromov-Witten
theory, often conjectural), and symplectic geometry (e.g.  \cite{lzz}).

Ekedahl, Lando, Shapiro and Vainshtein have announced a remarkable
formula~(\cite{elsv}~Theorem~1.1; Theorem \ref{biggie} below) linking
Hurwitz numbers to Hodge integrals in a particularly elegant way.

We prove Theorem \ref{biggie} using virtual localization on the moduli
space of stable maps, developed in \cite{gp}.  In the simplest case, no
complications arise, and Theorem \ref{biggie} comes out immediately;
Fantechi and Pandharipande proved this case independently
(\cite{fp}~Theorem 2), and their approach inspired ours.

We have chosen to present this proof because the formula of Ekedahl et al
is very powerful (see Sections \ref{mushroom} and \ref{celery} for
applications), and the program they propose seems potentially very
difficult to complete (e.g. \cite{elsv} Prop. 2.2, where they require
a compactification of the space of branched covers, with specified
branching at infinity, which is a bundle over $\cmbar_{g,n}$, such
that the branch map extends to the compactification).

In Section \ref{relative}, we show that the proof would be much
simpler if there were a moduli space for ``relative maps'' in the
algebraic category (with a good two-term obstruction theory, virtual
fundamental class, and hence virtual localization formula).  A space
with some of these qualities already exists in the symplectic category
(see \cite{lr} Section 7 and \cite{ip} for discussion).  In the algebraic case,
not much is known, although Gathmann has obtained striking results in
genus 0 (\cite{g}).

\epoint{Acknowledgements}  We are 
grateful to Rahul Pandharipande, David M. Jackson, and Michael Shapiro
for helpful conversations.

\section{Definitions and statement}

\point \label{avocado} \remind{avocado}
Throughout, we work over $\com$, and we use the following notation.
Fix a genus $g$, a degree $d$, and a partition $(\al_1,\ldots,\al_m)$
of $d$ with $m$ parts.  Let $b=2d+2g-2$, the ``expected number of
branch points of a degree $d$ genus $g$ cover of $\proj^1$'' by
the Riemann-Hurwitz formula.
We will identify $\Sym^b \proj^1$ with $\proj^b$ throughout.  Let
$r = d+m+2(g-1)$,
so a branched cover of $\proj^1$, with monodromy above
$\infty$ given by $\al$, and $r$ other specified simple branch points
(and no other branching) has genus $g$.
Let $k = \sum_i (\al_i-1)$, so $r=b-k$.
Let $H^g_{\al}$ be the number of such branched covers that are connected. 
(We do not take the points over $\infty$ to be labelled.)

\tpoint{Theorem (Ekedahl-Lando-Shapiro-Vainshtein, \cite{elsv}~Theorem 1.1)}
\label{biggie}
 \remind{biggie} {\em 
Suppose $g$, $m$ are integers ($g \geq 0$, $m \geq 1$) such that $2g-2+m>0$ (i.e.
the functor $\cmbar_{g,m}$ is represented by a Deligne-Mumford stack).  Then
\begin{equation*}
H^g_\al = \frac {r!} { \# \Aut(\al)}
 \prod_{i=1}^m \frac {{\al_i}^{\al_i}} {\al_i!} 
\int_{\cmbar_{g,m}}      \frac { 1-\la_1 + \dots \pm \la_g} {\prod (1-\al_i \psi_i)}
\end{equation*}
where $\la_i=c_i(\E)$ ($\E$ is the Hodge bundle).}

Fantechi and Pandharipande's argument applies in the case where there is no
ramification above $\infty$, i.e.  $\al = (1^d)$.

The reader may check that a variation of our method also shows that 
$$
H^0_{\al_1} = r! \frac {d^{d-2}} {d!}, \; \; \; H^0_{\al_1,\al_2} = \frac {r!}
{\# \Aut(\al_1, \al_2)} \cdot \frac {\al_1^{\al_1}} {\al_1!} \cdot
\frac {\al_2^{\al_2}} {\al_2!} \cdot d^{d-1}.
$$
As these formulas are known by other means (\cite{d} for the first, \cite{a} for the second,
\cite{gj0} for both), we omit the proof.

\epoint{Application:  Hurwitz numbers to Hodge integrals}
\label{mushroom} \remind{mushroom}
(i) Theorem \ref{biggie} provides a way of computing all Hodge integrals
as follows.  Define $$\langle \al_1, \dots, \al_m \rangle := \int_{\cmbar_{g,m}}
\frac { 1-\la_1 + \dots \pm \la_g} {\prod (1-\al_i \psi_i)},$$ a
symmetric polynomial in the $\al_i$ of degree $3g-3+m$ whose
coefficients are of the form $\int_{\cmbar_{g,m}} \psi_1^{d_1} \dots
\psi_m^{d_m} \la_k$.  It is straightforward to recover the
coefficients of a symmetric polynomial in $m$ variables of known
degree from a finite number of values, and $\langle \al_1, \dots, \al_m \rangle$ can
easily be computed (as Hurwitz numbers are combinatorial objects that
are easily computable, see Section \ref{pea}).  Once these integrals are
known, all remaining Hodge integrals (i.e. with more $\la$-classes) can be
computed in the usual way (\cite{m}).  The only other methods known to us are
Kontsevich's theorem, formerly Witten's conjecture,
\cite{ko}, which has no known algebraic proof,
and methods of Faber and Pandharipande (making clever use of 
virtual localization, \cite{p}).  These methods of computation
are in keeping with an extension of Mumford's philosophy, which is that 
much of the cohomology of $\cmbar_{g,n}$ is essentially combinatorial.

(ii) Combinatorially straightforward relations among Hurwitz numbers
(e.g. ``cut-and-join'', see \cite{gj0} Section 2) yield nontrivial new
identities among Hodge integrals.

\epoint{Application:  Hodge integrals to Hurwitz numbers} \label{celery}
\remind{celery}
There has been much work on the structure of the Hurwitz numbers,
including various predictions from physics.  Theorem \ref{biggie} is
the key step in a machine to verify these structures and predictions, 
see \cite{gjv}.

\section{Background:  Maps of curves to curves}
\point \label{chickpea} \remind{chickpea}
Following \cite{g01pn} Section 4.2, define a {\em special locus} of a
map $f: X \rightarrow \proj^1$ (where $X$ is a nodal curve) as
a connected component of the locus in $X$ where $f$ is not \etale.
(Remark:  No result in this section requires the target to be $\proj^1$.)
Then a special locus is of one of the following forms: (i) a nonsingular
point of $X$ that is an $m$-fold branch point (i.e. analytically
locally the map looks like $x \rightarrow x^m$, $m>1$), (ii) a node of
$X$, where the two branches of the node are branch points of order
$m_1$, $m_2$, or (iii) one-dimensional, of arithmetic genus $g$,
attached to $s$ branches of the remainder of the curve that are
$c_j$-fold branch points ($1 \leq j \leq s$).  The form of the locus,
along with the numerical data, will be called the {\em type}.  (For
convenience, we will consider a point {\em not} in a special locus to
be of type (i) with $m=1$.)  We will use the fact that special loci 
of type (ii) are smoothable (\cite{p1} Section 2.2).  

\epoint{Ramification number}
\label{broccoli} \remind{broccoli} 
To each special locus, associate a {\em ramification number} as
follows:  (i) $m-1$, (ii) $m_1 +m_2$, (iii) $2g-2+ 2s + \sum_{j=1}^s(c_j-1)$.
(Warning: in case (i), this is one less than what is normally called
the ramification index; we apologize for any possible confusion.)  The
{\em total ramification} above a point of $\proj^1$ is the sum of the
ramification numbers of the special loci mapping to that point.
We will use the following two immediate facts: if the map is stable, then
the ramification number of each ``special locus'' is a positive
integer, and each special locus of type (iii) has ramification number
at least 2.

\epoint{Extended Riemann-Hurwitz formula}
\label{pumpkin} \remind{pumpkin}
There is an easy generalization of the Riemann-Hurwitz formula:
$$
2 p_a(X) - 2 = -2d + \sum r_i
$$
where $\sum r_i$ is the sum of the ramification numbers.
(The proof is straightforward.  For example, consider the 
complex $f^* \omega^1_{\proj^1} \rightarrow \omega^1_X$ as
in \cite{fp} Section 2.3,
and observe that its degree can be decomposed into contributions from 
each special locus.  Alternatively, it follows from the usual 
Riemann-Hurwitz formula and induction on the number of nodes.)

\epoint{Behavior of ramification number and type in families} \label{yam} \remind{yam} Ramification
number is preserved under deformations.  Specifically, consider a
pointed one-parameter family of maps (of nodal cures).  Suppose one map
in the family has a special locus $S$ with ramification number $r$.
Then the sum of the ramification numbers of the special loci in a
general map that specialize to $S$ is also $r$.  (This can be shown by
either considering the complex $f^* \omega^1_{\proj^1} \rightarrow
\omega^1_X$ in the family or by deformation theory.)

Next, suppose 
$$\begin{array} {cc}
 C & \rightarrow \proj^1 \\
 \downarrow & \\
 B &
\end{array}
$$ 
is a family of {\em stable} maps parametrized by a nonsingular curve $B$.

\tpoint{Lemma}  \label{radish} \remind{radish}
{\em Suppose there is a point $\infty$ of $\proj^1$ where the total
ramification number of special loci mapping to $\infty$ is a constant
$k$ for all closed points of $B$.  Then the type of ramification above
$\infty$ is constant, i.e. the number of preimages of $\infty$ and their types are constant.}

For example, if the general fiber is nonsingular, i.e. only has
special loci of type (i), then that is true for all fibers.

{\em Proof.} Let $0$ be any point of $B$, and let $f: X \rightarrow
\proj^1$ be the map corresponding to $0$.  We will show that the type
of ramification above $\infty$ for $f$ is the same as for the general
point of $B$.

First reduce to the case where the general map has no contracted
components.  (If the general map has a contracted component $E$, then consider
the complement of the closure of $E$ in the total general family.  Prove the result
there, and then show that the statement of Lemma \ref{radish} behaves well
with respect to gluing a contracted component.)

Similarly, next reduce to the case where general map is nonsingular.
(First show the result where the nodes that are in the closure of the
nodes in the generic curve are normalized, and then show that the
statement behaves well with respect to gluing a 2-section of
the family to form a node.)

Pull back to an \etale neighborhood of 0 to separate special loci of
general fiber (i.e. so they are preserved under monodromy), and also
the fibers over $\infty$ for the general map.

For convenience of notation, restrict attention to one special locus
$E$ of $f$.  Assume first that $E$ is of type (iii), so $\dim E = 1$.
Let $g_E$ be the arithmetic genus of $E$.  Suppose that $r$ preimages
of $\infty$ of the general fiber (of type (i) by reductions) meet $E$
in the limit, and that these have ramification numbers $b_1$, \dots,
$b_r$.  Let $s$ be the number of other branches of $X$ meeting $E$, and $c_1$, 
\dots, $c_s$ the ramification numbers of the branches (as in
Section \ref{chickpea}).

The ramification number of $E$ is $(2 g_E-2) + 2s + \sum_{j=1}^s (c_j-1)$.
The total ramification number of the special loci specializing to $E$ is
$\sum_{i=1}^r (b_i-1)$.  
Also, 
$$\sum_{i=1}^r b_i = \sum_{j=1}^s c_j.$$
Hence by conservation
of ramification number, 
$$
(2 g_E-2+s) + r = 0.
$$
But $r>0$, 
and by the stability condition for $f$, $2g_E-2+s>0$, so we 
have a contradiction.

If $\dim E= 0$ is 0 (i.e. $E$ is of type (i) or (ii)), then
essentially the same algebra works (with the substitution ``$g_E=0$'',
resulting in $r+s-2=0$, from which $r=s=1$, from which the type is
constant).
\epf

A similar argument shows:
\tpoint{Lemma} \label{cucumber} \remind{cucumber}
{\em Suppose $E$ is a special locus in a specific fiber, and only one special locus $E'$ in the general fiber meets it.
Then the types of $E$ and $E'$ are the same.}

\epoint{The Fantechi-Pandharipande branch morphism}  
For any map $f$ from a nodal curve to a nonsingular curve, the ramification number defines a
divisor on the target: $\sum_L r_L f(L)$, where $L$ runs through the
special loci, and $r_L$ is the ramification number.  This induces a
set-theoretic map $\bbr: \cmbar_g(\proj^1,d) 
\rightarrow \Sym^b \proj^1 \cong \proj^b$.   
In \cite{fp}, this was shown to be a morphism.

Let $p$ be the point of $\Sym^b \proj^1 \cong \proj^b$
corresponding to $k(\infty) + (b-k)(0)$, let $L_\infty \subset
\proj^b$ be the linear
space corresponding to points of the form $k (\infty)+D$ (where $D$ is a divisor of degree $r=b-k$), and let
$\iota: L_\infty \rightarrow \proj^b$ be the inclusion.

Define $M$ as the stack-theoretic pullback $\bbr^{-1}L_\infty$.  It carries a virtual fundamental class
$[M]^{\vir} = \iota^! [ \cmbar_g(\proj^1,d) ]^{\vir}$ of dimension $r=b-k$
(i.e. simply intersect the class $[\cmbar_g(\proj^1,d)]^{\vir}$
with the
codimension $k$ operational Chow class $\bbr^*[L_\infty]$; the result
is supported on $\bbr^{-1}L_\infty$).
Denote the {\em restricted branch map} by $\rbr:  M \rightarrow L_\infty$. 
By abuse of notation, 
we denote the top horizontal arrow in the following diagram by $\iota$ as well.  
$$
\begin{array}{ccc}
M & \rightarrow & \cmbar_g(\proj^1,d) \\
\rbr \downarrow & & \downarrow \bbr \\
L_\infty & \stackrel \iota {\rightarrow} & \proj^b 
\end{array}
$$
By the projection formula,
\begin{equation}
\label{leek}
\iota_* ( \rbr^*[p] \cap [M]^{\vir}) = \bbr^*[p] \cap [ \cmbar_g(\proj^1,d) ]^{\vir}.\end{equation}
\remind{leek}

Define $M^\al$ as the union of irreducible components of $M$ whose
general members correspond to maps from irreducible curves, with
ramification above $\infty$ corresponding to $\al$ with the reduced substack
structure.  (It is not hard
to show that $M^\al$ is irreducible, by the same group-theoretic
methods as the classical proof that the Hurwitz scheme is irreducible.  None of
our arguments use this fact, so we will not give the details of the
proof.  Still, for convenience, we will assume irreducibility in our
language.) 

\point \label{kidneybean} \remind{kidneybean}
Note that $M=\bbr^{-1} L_\infty$ contains $M^\al$ with some multiplicity $m_\al$,
as $M^\al$ is of the expected dimension $r$. 
The Hurwitz number $H^g_\al$ is given by
$$\int_{M^\al} \rbr^* [p].$$  
(The proof of \cite{fp} Proposition 2 carries over without change in this case, as does 
the the argument of \cite{g1} Section 3.)
This is $1/ m_\al$ times the cap product of $\rbr^* [p]$ with 
the part of the 
class of $[M]^{\vir}$  supported on $M^\al$.

\tpoint{Lemma}  \label{lettuce} \remind{lettuce}
{\em $m_\al =  k! \prod \left( \frac {\al_i^{\al_i-1}} {\al_i!} \right).$}

\point \label{pea} \remind{pea} 
In the proof, we will use the combinatorial interpretation of Hurwitz
numbers: $H^g_\al$ is $1/d!$ times the number of ordered $r$-tuples
$(\tau_1, \dots, \tau_r)$ of transpositions generating $S_d$, whose
product has cycle structure $\al$.

{\em Proof.} Fix $r$ general points $p_1$, \dots, $p_r$ of $\proj^1$.
Let $L \subset \proj^b$ be the linear space corresponding to divisors
of the form $p_1 + \dots + p_r + D$ (where $\deg D = k$). By the
Kleiman-Bertini theorem, $(\bbr|_{M_\al})^{-1} L$ consists of $H^g_\al$
reduced points.  

Now $L_\infty \subset \Sym \proj^b$ can be interpreted as a ({\em
real} one-parameter) degeneration of the linear space 
corresponding to divisors of the form $D'+ \sum_{i=1}^k
q_i$, where $q_1$, \dots $q_k$ are fixed generally chosen points of
$\proj^1$ and $D'$ is any degree $r$ divisor on $\proj^1$. 

Choose branch cuts to the points $p_1$, \dots, $p_r$, $q_1$, \dots,
$q_k$, $\infty$ from some other point of $\proj^1$.  Choose a real
one-parameter path connecting $q_1$, \dots, $q_k$, $\infty$ (in that
order), not meeting the branch cuts (see the dashed line in Figure
\ref{degenfig}).  Degenerate the points $q_i$ to $\infty$ along this
path one at a time (so the family parametrizing this degeneration is
reducible).  If $\si_1$, \dots, $\si_k$, $\si_\infty$ are the
monodromies around the points $q_1$, \dots, $q_k$,
$\infty$ for a certain cover, then the monodromy around $\infty$ after
the branch points $q_i$, \dots, $q_k$ have been degenerated to
$\infty$ (along the path) is $\si_i \dots \si_k \si_\infty$.

At a general point of the family parametrizing this real degeneration  (before any of the
points $q_i$ have specialized, i.e. the $q_i$ are fixed general
points), $\bbr^{-1} (L \cap L_\infty)$ is a finite number of reduced
points.  This number is the Hurwitz number $H^g_{(1^d)}$ (\cite{fp}
Prop. 2), i.e.  $1/d!$ times the number of choices of $b=r+k$
transpositions $\tau_1$, \dots, $\tau_r$, $\si_1$,
\dots, $\si_k$ in $S_d$ such that $\tau_1 \dots \tau_r \si_1 \dots
\si_k$ is the identity and $\tau_1$, \dots, $\tau_r$, $\si_1$, \dots,
$\si_k$ generate $S_d$.

\begin{figure}
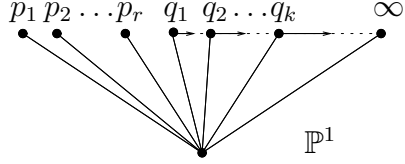

\begin{center}

	   \setlength{\unitlength}{.1\baseunit}
	    \input degenfig.tex 
\end{center}
\caption{Degenerating the points $q_1$, \dots, $q_k$ to $\infty$ one by one, along a real path \remind{degenfig}}
\label{degenfig}
\end{figure}

As we specialize the $k$ branch points $q_1$, \dots, $q_k$ 
to $\infty$ one at a time, some
of these points tend to points of $M^\al$; these are the points for
which $\tau_1$, \dots, $\tau_r$ generate $S_d$, and their product has
cycle structure $\al$.  The multiplicity $m_\al$ is the number of
these points that go to each point of $M^\al$.  This is the number of
choices of $k$ transpositions $\si_1$, \dots, $\si_k$ whose product is
a {\em given} permutation $\xi$ with cycle structure $\al$.  (Note that this 
number is independent of the choice of $\xi$; hence the multiplicity is
independent of choice of component of $M^\al$.)

If $k=\sum(\al_i-1)$ transpositions $\si_1$, \dots, $\si_k$ multiply
to a permutation $\xi=(a_{1,1} \dots a_{1,\al_1}) \dots (a_{m,1} \dots
a_{m, \al_m})$ (where $\{ a_{1,1}, \dots, a_{m,\al_m} \} = \{ 1,
\dots, d \}$), then for $1 \leq i \leq m$, $\al_i-1$ of the
transpositions must be of the form $(a_{i,j} a_{i,k})$.  (Reason:
A choice of $k+1$ points $q_1$, \dots, $q_k$, $\infty$, of $\proj^1$
and the
data $\si_1$, \dots, $\si_k$, $\xi$ defines a degree $d$ branched cover
of $\proj^1$, simply branched above $q_j$ and with ramification type
$\al$ above $\infty$.  By the Riemann-Hurwitz formula, 
the arithmetic genus of
this cover is $1-m$; as the pre-image of $\infty$ contains $m$ smooth
points, the cover has at most $m$ components.  Hence the cover has
precisely $m$ components, each of genus 0.  The $i$th component is
simply branched at $\al_i-1$ of the points $\{ q_1, \dots, q_k \}$
away from $\infty$.)

The number of ways of factoring an $\al_i$-cycle into $\al_i-1$
transpositions is $\al_i^{\al_i-2}$ (straightforward; or see \cite{d}
or \cite{gj0} Theorem 1.1).  Hence $m_\al$ is the number of ways of partitioning
the $k$ points $q_1$, \dots, $q_k$ into subsets of size $\al_1-1$, \dots, $\al_m-1$,
times the number of ways of factoring the $\al_i$-cycles:
$$
m_\al = \binom {k} {\al_1 -1, \dots, \al_m-1} \prod \al_i^{\al_i-2} = 
k! \prod \left( \frac {\al_i^{\al_i-1}} {\al_i!} \right).
$$
\epf

\section{Virtual localization}
\bpoint{Virtual localization preliminaries}
We evaluate the integral using virtual localization (\cite{gp}).  The
standard action of $\com^*$ 
on $\proj^1$ (so that the action on the tangent space at $\infty$ has weight 1) induces a natural $\com^*$-action 
on $\cmbar_g(\proj^1,d)$, and the branch morphism $\bbr$  is
equivariant with respect to the induced torus action on $\Sym^b \proj^1
\cong \proj^b$.  As a result, we can regard $\rbr^*[p]$ as an
equivariant Chow cohomology class in $A^r_{\com^*} M$.  Let $\{ F_l
\}_{l \in L}$ be the set of components of the fixed locus of
the torus action on $\cmbar_g(\proj^1,d)$, where $L$ is some index set.
(Note that the connected components of the fixed locus are also irreducible.)

Define $F_0$ to be the component of the fixed locus whose general point
parametrizes a
stable map with a single genus $g$ component contracted over 0, and
$m$ rational tails mapping with degree $\al_i$ ($1 \leq i \leq m$) to
$\proj^1$, totally ramified above 0 and $\infty$.  $F_{0}$ is naturally
isomorphic to a quotient of $\cmbar _{g,m}$ by a finite group.  See
\cite{k} or \cite{gp} for a discussion of the structure of the
fixed locus of the $\com ^{*}$ action on $\cmbar _{g}(\proj ^{1},d).$

By the virtual localization formula, we can explicitly write down classes $\mu_l \in A_*^{\com^*}(F_l)_{(1/t)}$
such that $$
\sum_l i_*(\mu_l) = [M]^{\vir}$$  in $A_*^{\com^*}(M)$.
Here, and elsewhere, $i$ is the natural inclusion.  It is important
to note that the $\mu_l$ are uniquely determined by this equation.
This follows from the Localization Theorem 1 of \cite{eg} (extended
to Deligne-Mumford stacks by \cite{kresch}), which says that pushforward gives
an isomorphism between the localized Chow group of the fixed locus and
that of the whole space.

In order to pick out the contribution to this integral from a 
single component $F_0$, we introduce more refined classes.  We
denote the irreducible components of $M$ by $M_n$, and arbitrarily 
choose a representation 
$$
[M]^{\vir} = \sum_n i_* \Ga_n
$$
where $\Ga_n \in A_*^{\com^*}( M_n)$. 
For a general
component, we can say little about these classes, but for our distinguished irreducible component $M^\al$
the corresponding $\Ga_\al$ is necessarily $m_\al [M^\al]$.  (Note that 
$M^\al$ has the expected dimension, so the Chow group in that 
dimension is generated by the fundamental class).

Next, we localize each of the $\Ga_n$.  Define $\eta_{l,n}$ in $A_*^{\com^*}(F_l)_{(1/t)}$ by
\begin{equation}
\label{rutabaga}
\sum_l i_* \eta_{l,n}= \Ga_n
\end{equation}
\remind{rutabaga}
Once again (by \cite{eg}, \cite{kresch}), the $\eta_{l,n}$ are uniquely defined; this will be used in Lemma \ref{potato}. 
Also, $\sum_n \eta_{l,n} = \mu_l$ (as the $\mu_l$ are uniquely determined).

\tpoint{Lemma} \label{squash}  \remind{squash}  {\em
The equivariant class $\rbr^*[p]$ restricts to 0 on any component
of the fixed locus whose general map has total ramification number greater
than $k$ above $\infty$.}

{\em Proof.}
Restricting the branch morphism to such a component, we see that it
gives a constant morphism to a point in $\proj ^{b}$ other than $p$.
Consequently, the pull-back of the class $p$ must vanish.
\epf

\tpoint{Lemma} \label{melanzana}  \remind{melanzana}  {\em 
$\int_{\Ga_n} \rbr^*[p] = 0$ for any irreducible component $M_n$ whose
general point corresponds to a map which has a contracted 
component away from $\infty$.}

{\em Proof.} A general cycle $\gamma \in L_\infty$ representing $p$ is
the sum of $r$ distinct points plus the point $\infty$ exactly $k$
times.  However, a contracted
component always gives a multiple component of the branch divisor,
(Section \ref{broccoli}) so the image of $M_{n}$ cannot meet a general 
point.
\epf

\tpoint{Lemma} {\em $\eta_{l,n} = 0$ if $F_l \cap M_n = \emptyset$.}
\label{potato} \remind{potato}

{\em Proof.} 
Since $\Ga_n$ is an element of $A_*^{\com^*}(M_n)$, there exist
classes $\tilde{\eta}_{l,n}$ in the localized equivariant Chow groups
of the fixed loci of $M_n$ satisfying equation (\ref{rutabaga}).
Pushing these forward to the fixed loci of $M$ gives classes in the Chow groups of the $F_l$
satisfying the same equation.  By uniqueness, these must be the $\eta_{l,n}$.  By this
construction, it follows that they can only be non-zero if $F_l$ meets $M_n$. \epf

\tpoint{Lemma} {\em No irreducible component of $M$ can meet two distinct components of the 
fixed locus with total ramification number exactly $k$ above $\infty$.}
\label{cabbage}
\remind{cabbage}

{\em Proof.} To each map $f: X \rightarrow \proj^1$ with total ramification
number exactly $k$ above $\infty$, associate a graph as follows.  The
connected components of the preimage of $\infty$ correspond to red
vertices; they are labelled with their type.  The connected components of $Y=
\overline{X \setminus f^{-1}( \infty )}$ (where the closure is taken in
$X$) correspond to green vertices; they are labelled with their
arithmetic genus.  Points of $Y \cap f^{-1}( \infty)$ correspond to
edges connecting the corresponding red and green points; they are
labelled with the ramification number of $Y \rightarrow \proj^1$ at
that point.  Observe that this associated graph is constant in
connected families where the total ramification over $\infty$
is constant, essentially by Lemma \ref{radish}.

If an irreducible component $M'$ of $M$ meets a component of the fixed
locus with total ramification number exactly $k$ above $\infty$, then the
general map in $M'$ has total ramification $k$ above $\infty$.
(Reason: the total ramification is at most $k$ as it specializes to a
map with total ramification exactly $k$; and the total ramification is
at least $k$ as it is a component of $M$.)  There is only one
component of the fixed locus that has the same associated graph as the
general point in $M'$, proving the result. \epf

\tpoint{Lemma} \label{arugula} \remind{arugula}  {\em 
The map parametrized by a general point of any 
irreducible component of $M$ other
than $M^\al$ which meets $F_0$ must have a contracted component not mapping
to $\infty$.} 

{\em Proof.} Let $M'$ be an irreducible component of $M$ other than
$M^\al$.  As in the proof of Lemma \ref{cabbage}, a general map $f:  X \rightarrow \proj^1$ of
$M$ has total ramification exactly $k$ above $\infty$.  By Lemma
\ref{radish}, we know the type of the special loci above $\infty$:
they are nonsingular points of the source curve, and the ramification
numbers are given by $\al_1, \dots, \al_m$.

As $M' \neq M$, $X$ is singular.  If $f$ has a special locus of type
(iii), then we are done.  Otherwise, $f$ has only special loci of type
(ii), and none of these map to $\infty$.  But then these type (ii)
special loci can be smoothed while staying in $M$ (Section
\ref{broccoli}), contradicting the assumption that $f$ is a general
map in a component of $M$.
\epf

\tpoint{Proposition} 
\label{onion} \remind{onion} $$
m_\al \int_{M^{\al}} \rbr^*[p] = \int_{F_0} \rbr^* [p] \cap \mu_0.$$

It is the class  $\mu_0$ that the Virtual Localization Theorem 
of \cite{gp}
allows us to calculate explicitly.  Thus this proposition is the main ingredient in giving us
an explicit formula for the integral we want to compute.

{\em Proof.}
Now $\Ga_\al = m_\al [M^\al]$, so by definition of $\eta_{l,\al}$,
$$
m_\al [M^\al] = \sum_l i_* \eta_{l,\al}.$$

By Lemma \ref{melanzana}, $M^\al$ meets only one component of the fixed locus which
has total ramification number $k$, $F_0$.  Along with Lemmas \ref{squash} 
and \ref{potato},
this implies that
$$
m_\al \int_{M^\al} \rbr^*[p] = \int_{F_0} \rbr^*[p] \cap \eta_{0,\al}.
$$
In other words, the only component of the fixed locus which
contributes to this integral is $F_0$.  Since $\mu_0 = \sum_n
\eta_{0,n}$, the proposition will follow if we can show that
$$\int_{F_0} \rbr^* [p] \cap \eta_{0,n} = 0$$ for $n \neq \al$, i.e. that
no other irreducible component of $M$ contributes to
the localization term coming from $F_0$.  

If $F_0 \cap M_n = \emptyset$, this is true by Lemma \ref{potato}.
Otherwise, by Lemma \ref{arugula}, the general map in $M_n$
has a contracted  component, so by Lemma \ref{melanzana} 
$\int_{\Gamma_n} \rbr^*[p] = 0$.  By equation (\ref{rutabaga}),
$$
\sum_l \int_{F_l} \rbr^*[p] \cap \eta_{l,n} = 0.$$
If $F_l$ generically corresponds to maps that have total ramification
number greater than $k$ above $\infty$, then $\rbr^*[p] \cap \eta_{l,n} = 0$
by Lemma \ref{melanzana} .  If $l \neq 0$ and $F_l$ generically corresponds
to maps that have total ramification number $k$ above $\infty$,
then $\rbr^*[p] \cap \eta_{l,n} = 0$ by Lemma \ref{cabbage}, as $M_n$ meets
$F_0$.  Hence $\int_{F_0} \rbr^*[p] \cap \eta_{0,n} = 0$ as desired.
\epf

\bpoint{Proof of Theorem \ref{biggie}}
\label{proof} \remind{proof}
All that is left is to explicitly write down the right hand side of
Proposition \ref{onion}.  By equation (\ref{leek}), this integral
can be interpreted as the contribution of $F_0$ to the integral of
$\bbr^*[p]$ against the virtual fundamental class of $\cmbar_g(\proj^1,d)$, divided by
$m_\al$.   Since this means we are trying to compute an equivariant
integral over the entire space of maps to $\proj^1$, we are in exactly
the situation discussed in \cite{gp}.  Let $\gamma$ be the natural
morphism from $\cmbar_{g,m}$ to $F_0$.  The degree of $\gamma$ is $\#
\Aut(\alpha)\prod \alpha_i$.  The pullback under $\gamma$ of the
inverse euler class of the virtual normal bundle is computed to be
$$c(E^\vee) \left( \prod \frac{1}{1-\alpha_i\psi_i} \cdot
\frac{(-1)^{\alpha_i}\alpha_i^{2\alpha_i}}{(\alpha_i!)^2} \right).$$
The class $\rbr^*[p]$ is easy to evaluate.  Since $\rbr$ is constant when 
restricted to $F_0$, this class is pure weight, and is given by the 
product of the weights of the $\com^*$ action on $T_p \proj^b$.  These 
weights are given by the non-zero integers from $-(b-k)$ to $k$ inclusive.
The integral over $F_0$ is just the integral over $\cmbar_{g,m}$ divided 
by the degree of $\gamma$.  We conclude that 
$$m_\al \int_{[M^\alpha]} br^*[p] =
\frac{k!(b-k)!}{\# \Aut(\alpha) \prod \alpha_i} \cdot \prod 
\frac{\alpha_i^{2\alpha_i}}{(\alpha_i!)^2} \cdot \int_{\cmbar_{g,m}} 
\frac{c(\E^\vee)}{\prod (1-\alpha_i \psi_i)}.
$$
Dividing by $m_\alpha$ (calculated 
in Lemma \ref{lettuce}) yields the desired formula.
\epf

\section{A case for an algebraic definition of a 
space of ``relative stable maps''}
\label{relative}
\remind{relative}

A space of ``relative stable maps'' has been defined in the symplectic category
(see \cite{lr} and \cite{ip}), but hasn't yet been properly defined in the algebraic
category (with the exception of Gathmann's work in genus 0, \cite{g}).

The proof of Theorem \ref{biggie} would become quite short were such a
space $\cm$ to exist with expected
properties, namely the following.  Fix $d$, $g$, $\al$, $m$, $k$, $r$ as
before (see Section \ref{avocado}).
\begin{enumerate}
\item $\cm$ is a proper Deligne-Mumford stack, which contains as an open substack $U$ the
locally closed substack of $\cmbar_g(\proj^1,d)$ corresponding to maps to $\proj^1$ where
the pre-image of $\infty$ consists of $m$ smooth points appearing with multiplicity
$\al_1$, \dots, $\al_m$.
\item There is a Fantechi-Pandharipande branch map $\bbr:  \cm \rightarrow
\Sym^b \proj^1$.  The image will be contained in $L_\infty$, so we may consider
the induced map $\rbr$ to $L_\infty \cong \Sym^r \proj^1$.  Under this
map, the set-theoretic fiber of $k(\infty) +r(0)$ is precisely $F_{0}$.
\item There is a $\com^{*}$-equivariant perfect obstruction theory on 
$\cm $ which when
restricted to $U$ is given (relatively over ${\mathfrak{M}}_g$) by 
$R\pi _{*}(f^{*}(T\proj ^{1} \otimes \oh (- \infty)))$, where
$\pi$ is the structure morphism from the universal curve to $\cm$.
\end{enumerate}
With these axioms, the proof would require only Section \ref{proof}.

All of these requirements are reasonable.  However, as a warning, note that
the proof of Proposition \ref{onion} used special properties of the class 
$\rbr^* [p]$ (Lemmas \ref{squash}--\ref{arugula}).

One might expect this space to be a combination of Kontsevich's space
$\cmbar_g(\proj^1,d)$ and the space of twisted maps introduced by 
Abramovich and Vistoli (see \cite{av} Section 3).

} 

\end{document}

%% file: degenfig.tex
\begingroup\makeatletter\ifx\SetFigFont\undefined%
\gdef\SetFigFont#1#2#3#4#5{%
  \reset@font\fontsize{#1}{#2pt}%
  \fontfamily{#3}\fontseries{#4}\fontshape{#5}%
  \selectfont}%
\fi\endgroup%
{\renewcommand{\dashlinestretch}{30}
\begin{picture}(6944,2654)(0,-10)
\put(225,2183){\blacken\ellipse{150}{150}}
\put(225,2183){\ellipse{150}{150}}
\put(825,2183){\blacken\ellipse{150}{150}}
\put(825,2183){\ellipse{150}{150}}
\put(2025,2183){\blacken\ellipse{150}{150}}
\put(2025,2183){\ellipse{150}{150}}
\put(3525,2183){\blacken\ellipse{150}{150}}
\put(3525,2183){\ellipse{150}{150}}
\put(4725,2183){\blacken\ellipse{150}{150}}
\put(4725,2183){\ellipse{150}{150}}
\put(6525,2183){\blacken\ellipse{150}{150}}
\put(6525,2183){\ellipse{150}{150}}
\put(3375,83){\blacken\ellipse{150}{150}}
\put(3375,83){\ellipse{150}{150}}
\put(2872,2205){\blacken\ellipse{150}{150}}
\put(2872,2205){\ellipse{150}{150}}
\path(3375,83)(225,2183)
\path(3375,83)(225,2183)
\path(3375,83)(825,2183)
\path(3375,83)(825,2183)
\path(2025,2183)(3375,83)
\path(2025,2183)(3375,83)
\path(3525,2183)(3375,83)
\path(3525,2183)(3375,83)
\path(4725,2183)(3375,83)
\path(4725,2183)(3375,83)
\path(3375,83)(6525,2183)
\path(3375,83)(6525,2183)
\path(4800,2183)(5625,2183)
\path(4800,2183)(5625,2183)
\path(5505.000,2153.000)(5625.000,2183.000)(5505.000,2213.000)
\path(3525,2183)(4125,2183)
\path(3525,2183)(4125,2183)
\path(4005.000,2153.000)(4125.000,2183.000)(4005.000,2213.000)
\path(3375,83)(2850,2183)
\path(3375,83)(2850,2183)
\path(2850,2183)(3225,2183)
\path(2850,2183)(3225,2183)
\path(3105.000,2153.000)(3225.000,2183.000)(3105.000,2213.000)
\dashline{60.000}(2850,2183)(6525,2183)
\put(1875,2483){\makebox(0,0)[lb]{\smash{{{\SetFigFont{12}{14.4}{\rmdefault}{\mddefault}{\updefault}$p_r$}}}}}
\put(6375,2483){\makebox(0,0)[lb]{\smash{{{\SetFigFont{12}{14.4}{\rmdefault}{\mddefault}{\updefault}$\infty$}}}}}
\put(0,2483){\makebox(0,0)[lb]{\smash{{{\SetFigFont{12}{14.4}{\rmdefault}{\mddefault}{\updefault}$p_1$}}}}}
\put(600,2483){\makebox(0,0)[lb]{\smash{{{\SetFigFont{12}{14.4}{\rmdefault}{\mddefault}{\updefault}$p_2$}}}}}
\put(1200,2483){\makebox(0,0)[lb]{\smash{{{\SetFigFont{12}{14.4}{\rmdefault}{\mddefault}{\updefault}$\dots$}}}}}
\put(3375,2483){\makebox(0,0)[lb]{\smash{{{\SetFigFont{12}{14.4}{\rmdefault}{\mddefault}{\updefault}$q_2$}}}}}
\put(4575,2483){\makebox(0,0)[lb]{\smash{{{\SetFigFont{12}{14.4}{\rmdefault}{\mddefault}{\updefault}$q_k$}}}}}
\put(5175,83){\makebox(0,0)[lb]{\smash{{{\SetFigFont{12}{14.4}{\rmdefault}{\mddefault}{\updefault}$\proj^1$}}}}}
\put(2700,2483){\makebox(0,0)[lb]{\smash{{{\SetFigFont{12}{14.4}{\rmdefault}{\mddefault}{\updefault}$q_1$}}}}}
\put(3900,2483){\makebox(0,0)[lb]{\smash{{{\SetFigFont{12}{14.4}{\rmdefault}{\mddefault}{\updefault}$\dots$}}}}}
\end{picture}
}